\documentclass[12pt]{amsart}
\usepackage{amsfonts,amssymb}
\usepackage{mathrsfs}

\usepackage[all]{xy}
\usepackage{hypbmsec}
\usepackage{mathtext}

{}

\numberwithin{equation}{section}

\newtheorem{theorem}{Theorem}[section]
\newtheorem{theoremnad}{Nadel Vanishing Theorem}
\newtheorem{lemma}[theorem]{Lemma}
\newtheorem{prop}[theorem]{Proposition}
\newtheorem{corollary}[theorem]{Corollary}

 \theoremstyle{remark}
\newtheorem{remark}[theorem]{Remark}

\theoremstyle{definition}
\newtheorem{defin}[theorem]{Definition}
\newtheorem{example}[theorem]{Example}
\newenvironment{fulltext}{}{}

\begin{document}

\title
{On $G$-Fano threefolds}

\author{Yu.~G.~Prokhorov}
\thanks{This work is supported by the Russian Science Foundation under grant  14-50-00005}

\address{Steklov Mathematical Institute, Russian Academy of Sciences;
}
\email{prokhoro@gmail.com}

\maketitle


\begin{abstract}
We study Fano threefolds with~terminal
singularities admitting a ``minimal'' action of a finite group.
We prove that under certain additional assumptions such a variety does not contain
planes. We also obtain an upper bounds of the number of singular points of
certain Fano threefolds with~terminal factorial singularities.
\end{abstract}

 \begin{fulltext}


\section{Introduction}
\label{s1}

Let $X$~be an algebraic variety over a field $\Bbbk$ of characteristic zero
and let $G$~be some group. We say that ~$X$
is a $G$-\textit{variety}, if the group~$G$ acts on
$\overline X=X\otimes\bar\Bbbk$, where $\bar\Bbbk$~is the algebraic closure of~$\Bbbk$. 
Moreover, we assume that ~$X$, $G$ and~$\Bbbk$
satisfy one of the following conditions.

(a) \textit{Geometric case}: the field $\Bbbk$ algebraically closed,
the group~$G$ is finite and the action of~$G$ on~$X$ is defined by a homomorphism
$G\to\operatorname{Aut}_{\Bbbk}(X)$.

(b) \textit{Algebraic case}: $G$~is the Galois group of $\bar\Bbbk$ over
$\Bbbk$ acting on $\overline X=X\otimes\bar\Bbbk$ through the second factor. 
The action of~$G$ on~$X$ is trivial.

A $G$-variety $X$ is called a $G$-\textit{Fano variety}, if
the singularities of~$X$ are not worse than terminal Gorenstein, the
anticanonical
divisor~$-K_X$ is ample and the rank of the invariant part~$\operatorname{Cl}(X)^G$
of the Weil divisor class group~$\operatorname{Cl}(X)$ equals~$1$
(see~\cite{1}--\cite{3}). In the present paper we consider only
the three-dimensional case.

We say that a Fano threefold $X$ \textit{belongs to the main
series}, if its canonical divisor~$K_X$ generates the Picard group
$\operatorname{Pic}(X)$. $G$-Fano threefolds of non-main series were classified
by the author in the works~\cite{2},~\cite{4}.

Recall that \textit{the genus} of a Fano threefold $X$ is the number 
$\mathrm g(X):=\frac12(-K_X)^3+1$ (see~Definition~\ref{d2.1}).

We prove the following theorem.

\begin{theorem}
\label{t1.1}
Let $X$~be a $G$-Fano threefold of the main series
with~$\mathrm g(X)\ge6$. Then~$X$ does not contain any planes.
\end{theorem}

It turns out that the absence of planes on a Fano threefold
with~terminal singularities is very important for the classification. 
Indeed, for such a variety there exists a 
$\mathbb Q$-factorialization $\pi\colon X'\to \nobreak X$, where $X'$~is a variety 
with~terminal factorial singularities and numerically effective (nef) big
anticanonical divisor. If there are no planes on~$X$, then running 
the minimal model program to~$X'$ we stay in the same class of varieties 
(the category of terminal factorial varieties with~nef and big anticanonical divisor which does not contain planes; see~\cite{5}, \cite{6}). 
In some cases (especially for large values of genus) this allows 
to obtain a description of the original variety~$X$.
Applications of Theorem~\ref{t1.1} that use this construction will be discussed 
in the forthcoming paper.

Note that for small values of genus, $G$-Fano threefolds can contain planes.

\begin{example}
\label{e1.2}
\textit{The Burkhardt quartic} $X_4^{\mathrm b}$~is the subvariety
in~$\mathbb P^5$, defined by the equations $\sigma_1=\sigma_4=0$, where
$\sigma_i$~are elementary symmetric polynomials in $x_1,\dots,x_6$. This 
quartic was intensively studied earlier (see, e.g., \cite{7}). The singular
locus of~$X_4^{\mathrm b}$ consists of~45 ordinary double points.
The symmetric group~$\mathfrak S_6$ acts on~$X_4^{\mathrm b}$
by permutations of coordinates. Then the 
quotient variety~$\mathbb P^5/\mathfrak S_6$ is isomorphic to the weighted projective space
$\mathbb P(1,\dots,6)$ and the
quotient variety~$X_4^{\mathrm b}/\mathfrak S_6$ is isomorphic to the
subspace $\mathbb P(2,3,5,6)\subset\mathbb P(1,\dots,6)$.
Therefore,
$\operatorname{rk}\operatorname{Cl}(X_4^{\mathrm b})^{\mathfrak S_6}=1$,
and so $X_4^{\mathrm b}$ is a $\mathfrak S_6$-Fano threefold 
of the main series and genus~$3$. The quartic~$X_4^{\mathrm b}$ contains exactly ~40
planes~\cite{7}.
\end{example}

It is known also a lot of examples of Fano threefolds of large genus (which are not 
$G$-Fano) with~terminal singularities that
contain planes. However the author does not know any examples of
$G$-Fano threefolds of the main series of genus~$4$ and~$5$ containing planes 
(see~Corollary~\ref{c3.12}).

$\mathbb Q$-factorial (terminal) Fano threefolds of the main series
are always $G$-Fano with respect, for instance, to the trivial group. In this case,
for $\mathrm g(X)\ge8$, we obtain an upper bound for the number 
of singular points which is sharp for $\mathrm g(X)\ge9$.

\begin{theorem}
\label{t1.3}
Let $X$~be a $\mathbb Q$-factorial Fano threefold
of the main series with terminal singularities. Then, for
$\mathrm g(X)=9,10,12$, the number of singular points of~$X$ is 
at most $12-\mathrm g(X)$ and this bound is sharp. 
For $\mathrm g(X)=8$ the variety~$X$ has at most~$10$
singular points.
\end{theorem}

Moreover, we generalize the classical Fano--Iskovskikh ``double projection'' construction 
 (see~Theorem~\ref{t4.1}).

The paper is organized as follows: \S\,\ref{s2} is 
preliminary; in~\S\,\ref{s3} we prove Theorem~\ref{t1.1}; Theorem~\ref{t1.3} is 
deduced from more general Theorem~\ref{t4.1}
in~\S\,\ref{s4}; finally, \S\,\ref{s5} contains two auxiliary
results which are used in~the proof Theorem~\ref{t4.1}.

\smallskip

The author would like to thank the referee for constructive comments.

\section{Preliminaries}
\label{s2}
For certainty, considering $G$-varieties, we deal only with the
geometric case. The algebraic case is similar. Thus
the ground field~$\Bbbk$ is supposed to be algebraically closed
(of characteristic ~$0$). Sometimes we will assume also that 
$\Bbbk=\mathbb C$.

All the varieties considered in this paper have at worst terminal Gorenstein singularities. 
Every such a three-dimensional singularity $X\ni P$
is locally a hypersurface and has multiplicity~2. The Picard group of a 
variety with terminal singularities is embedded 
to~the Weil divisor class group so that the cokernel has no torsion elements~\cite[Lemma~5.1]{8}.

\begin{defin}
\label{d2.1}
Let $X$~be a Fano threefold with~terminal
singularities. Its \textit{genus} is the number
$\mathrm g(X):=-K_X^3/2+1$. 
\end{defin}

By the Riemann--Roch formula and the
 Kawamata--Viehweg vanishing theorem we have $\dim\lvert-K_X\rvert=\mathrm g(X)+1$.
In particular, $\mathrm g(X)$~is an integer. In the case of Fano threefolds
of the main series, the genus can take only the following values:
$\mathrm g(X)\in\{2,3,\dots,10,12\}$ (see~\cite{9} and~\cite{10}).

The Picard group $\operatorname{Pic}(X)$ and the Weil divisor class group
$\operatorname{Cl}(X)$ are finitely generated and torsion free
\cite{9}.

\begin{theorem}[{\rm\cite{11}, \cite{12}}]
\label{t2.2}
Let $X$~be a Fano threefold with~terminal singularities and
with~$\operatorname{Pic}(X)\simeq\mathbb Z\cdot K_X$. The following assertions 
hold:

{\rm(i)} the linear system $\lvert-K_X\rvert$ is base point free;

{\rm(ii)} if $\mathrm g(X)\ge4$, then $\lvert-K_X\rvert$ is very ample and defines
an embedding $X=X_{2g-2}\subset\mathbb P^{g+1}$;

{\rm(iii)} if $\mathrm g(X)\ge5$, then the image
$X=X_{2g-1}\subset\mathbb P^{g+1}$ is an intersection of quadrics.
\end{theorem}

\begin{theorem}[{\rm\cite{10}}]
\label{t2.3}
Let $X$~be a Fano threefold with~terminal
singularities. Then~$X$ is smoothable, i.\,e.~there exists a flat family
$\mathfrak f\colon\mathfrak X\to(\mathfrak D\ni0)$ over a disk
$(\mathfrak D\ni0)\subset\mathbb C$ such that $\mathfrak X_0\simeq X$ and a
general element~$\mathfrak X_s$, $s\in\mathfrak D$ is a nonsingular
Fano threefold. Moreover, there exist natural identifications
$\operatorname{Pic}(X)=\operatorname{Pic}(\mathfrak X_s)
=\operatorname{Pic}(\mathfrak X)$ so that 
$K_{\mathfrak X_s}=K_X$ (see~\cite[\S\,1]{13}).
\end{theorem}

Let $(X,B)$~be a log pair (a pair consisting of a normal variety~$X$
and an effective $\mathbb Q$-divisor $B=\sum_ib_iB_i$ on~$X$). Assume
that $K_X+B$~is a $\mathbb Q$-Cartier divisor. Let $f\colon\tilde X\to X$~ be a
log resolution $(X,B)$. Write
$$
K_{\tilde X}=f^*(K_X+B)+E,
$$
where $E=\sum_ie_i E_i$~is a $\mathbb Q$-divisor whose components are
proper transforms of the components of~$B$ and the exceptional divisors. By our
hypothesis $\sum_iE_i$ has only simple normal crossings. The pair
$(X,B)$ has \textit{log canonical (lc) singularities} if $e_i\ge-1$ for all~$i$.
A proper irreducible subvariety $Z\subset X$ is called
\textit{a center of log canonical singularities} $(X,B)$ if, for some
resolution~$f$, there exists a component~$E_i$ with~coefficient
$e_i\le-1$ dominating~$Z$. The union of all centers of
log canonical singularities is called \textit{the locus of log canonical
singularities} and denoted by $\operatorname{LCS}(X,B)$. 
Thus,
$$
\operatorname{LCS}(X,B)=\bigcup_{e_i\le-1}f(E_i).
$$
The sheaf
$$
\mathcal I(X,B):=f_*\mathscr O_{\tilde X}(\lceil E\rceil)
$$
is called \textit{the multiplier sheaf}. Since $B$ is effective, 
$\mathcal I(X,B)$ is an ideal sheaf. The corresponding subscheme
in~$X$ is called \textit{the scheme of log canonical singularities}. Its support
coincides with $\operatorname{LCS}(X,B)$. If the pair $(X,B)$ is lc, then
$\mathscr O_X/\mathcal I(X,B)$ has no nilpotents and so the scheme of
log canonical singularities is reduced (and coincides
with~$\operatorname{LCS}(X,B)$).

\begin{theoremnad}[{\rm (\cite[Theorem~9.4.17]{14})}]
Let $(X,B)$ be an lc pair, where the variety~$X$ is projective.
Let $D$~be a Cartier divisor on~$X$ such that the divisor $D-(K_X+B)$ is nef and big. Then
$$
H^q(\mathcal I(X,B)\otimes\mathscr O_X(D))=0 \quad
\forall\,q>0.
$$
\end{theoremnad}

\section{Planes}
\label{s3}

In this section we prove Theorem \ref{t1.1}.

First we introduce the notation. Let $X$~be a $G$-Fano threefold
of the main series with~terminal singularities. We
assume that $g=\mathrm g(X)\ge5$. Thus,
$\operatorname{Pic}(X)\simeq\mathbb Z\cdot K_X$ and the linear system
$\lvert-K_X\rvert$ defines an embedding $X=X_{2g-2}\subset\mathbb P^{g+1}$ so that 
its the image~is an intersection of quadrics (by Theorem \ref{t2.2}).

Assume that there exists a plane $\Pi_1\subset X$. Let
$O=\{\Pi_1,\dots,\Pi_{n}\}$~be its orbit with respect to the action of~$G$ and let
$D:=\sum_i\Pi_i$. Recall that $\operatorname{Cl}(X)^G\simeq
\operatorname{Pic}(X)^G\simeq\mathbb Z\cdot K_X$. Hence ~$D$ is a
Cartier divisor and for some integer~$a$ we can write
$D\sim-aK_X$. Comparing the degrees we obtain
\begin{equation}
\label{eq3.1}
n=(2g-2)a.
\end{equation}

It is clear that for any two distinct planes $\Pi_i,\Pi_j\in O$ their
intersection $\Pi_i\cap\Pi_j$ is either empty, a point or a line.

\begin{lemma}
\label{l3.1}
In the above notation, the number of planes passing through any point
$P\in X\setminus\operatorname{Sing}(X)$ (and contained in~$X$) is at most two. 
In particular, the divisor~$D$ has only simple normal
crossings in the nonsingular locus $X\setminus\operatorname{Sing}(X)$.
\end{lemma}

\begin{proof}
Let $P\in X$~be a nonsingular point and let $\Pi_1,\dots,\Pi_r\in O$~ be all the
planes passing through~$P$. Then these planes are contained
in the~projective tangent space
$\overline{T_{P,X}}\simeq\mathbb P^3$ to~$X$ at~$P$. Since
$X\subset\mathbb P^{g+1}$ is an intersection of quadrics (by
Theorem~\ref{t2.2}), the subvariety 
$\overline{T_{P,X}}\cap X$ is an intersection of
quadrics and so it cannot contain more than two planes.
\end{proof}

\begin{corollary}
\label{c3.2}
The pair $(X,D)$ has only log canonical singularities
in~$X\setminus\operatorname{Sing}(X)$. Moreover,
$X\setminus\operatorname{Sing}(X)$ does not contain any zero-dimensional log canonical
centers.
\end{corollary}

\begin{lemma}
\label{l3.3}
There are at most four planes passing trough a singular point $P\in X$
(and contained in~$X$).
\end{lemma}

\begin{proof} 
As in~Lemma \ref{l3.1}, all the planes $\Pi_1,\dots,\Pi_r$
passing through~$P$ are contained in~the set $\overline{T_{P,X}}\cap X$
which is an intersection of quadrics in~$\overline{T_{P,X}}$. Since
$P\in X$~is a hypersurface singularity, $\dim T_{P,X}=4$. Since
$\dim\overline{T_{P,X}}\cap X\le2$, we we obtain that $\overline{T_{P,X}}\cap X$
contains at most four planes.
\end{proof}

\begin{lemma}
\label{l3.4}
The pair $(X,D)$ is lc.
\end{lemma}

\begin{proof}
Assume the contrary. Then $(X,(1-\varepsilon)D)$ is not lc for
$0<\varepsilon\ll1$. According to Corollary~\ref{c3.2} the locus of
log canonical singularities $\operatorname{LCS}(X,(1-\varepsilon)D)$~ is a 
finite set of points (non-empty), it is contained in~ the singular locus of~$X$.
On the other hand, by Shokurov's connectedness theorem
(see~\cite[Theorem~17.4]{15}) this set is connected\footnote{
This argument works for $a=1$. For $a>1$ one can use the 
inversion of adjunction to a plane $\Pi_i$ (see \cite{16})
and the fact that $\overline{T_{P,X}}\cap X$ is an intersection of 
quadrics (see Theorem \ref{t2.2}).
}. 
Hence,
$\operatorname{LCS}(X,(1-\varepsilon)D)$~is a single point~$P$ which must
be $G$-invariant and singular for~$X$. Then \textit{all the components of}~$D$
pass through~$P$. This contradicts Lemma~\ref{l3.3} because the number of
components of~$D$ greater than~4. The lemma is proved.
\end{proof}

Since $D$~is a Cartier divisor on a variety with~terminal
singularities, $D$ is a Cohen--Macaulay scheme. Therefore, for
any component $\Pi_i\subset D$, the intersection
$\Pi_i\cap\operatorname{Supp}(D-\Pi_i)$ has pure dimension~$1$.
On the other hand, the scheme $\Pi_i\cap\operatorname{Supp}(D-\Pi_i)$ is reduced
in the~generic point by Lemma~\ref{l3.1} (and its components~are distinct lines).
Put $\Delta_i:=\Pi_i\cap\operatorname{Supp}(D-\Pi_i)$.

\begin{corollary}
\label{c3.5}
For any component $\Pi_i\subset D$, the divisor $\Delta_i\subset\Pi_i$ has
simple normal crossing.
\end{corollary}

\begin{proof}[{\rm follows from Shokurov's log canonical inversion of adjunction
 (see~\cite{16})}]
\end{proof}

It is clear that two-dimensional centers of log canonical singularities of the pair $(X,D)$
are planes~$\Pi_i$ and one-dimensional ones are those intersections
$\Pi_i\cap\Pi_j$ that are lines. Denote by $\mathscr P\subset X$~ the
set of all zero-dimensional centers of log canonical singularities of $(X,D)$.
According to Corollary~\ref{c3.2} we have
$\mathscr P\subset\operatorname{Sing}(X)$.

\begin{lemma}
\label{l3.6}
For any component $\Pi_i\subset D$ we have
$\mathscr P\cap\Pi_i=\operatorname{Sing}(\Delta_i)$.
\end{lemma}

\begin{proof}
Fix a point $P\in\Pi_i$. Let $H\subset X$~be a general hyperplane
section passing through~$P$.

Let $P\in\operatorname{Sing}(\Delta_i)$. If~$P$ is not a
log canonical center, then the pair $(X,D+\varepsilon H)$ is lc for
$0<\varepsilon\ll1$. In this case, by the inversion of adjunction~\cite{16} the pair
$(\Pi_i,\Delta_i +\varepsilon H|_{\Pi_i})$ is also lc which is impossible
because the multiplicity of $\Delta_i +\varepsilon H|_{\Pi_i}$ at~$P$
is greater than~$2$. The contradiction shows that 
$\mathscr P\supset\operatorname{Sing}(\Delta_i)$.

Conversely, let $P\in\mathscr P$. Again by the inversion of adjunction the pair
$(\Pi_i,\Delta_i+\varepsilon H|_{\Pi_i})$ is not lc for
$0<\varepsilon\ll1$. Therefore, the curve~$\Delta_i$ is singular at~$P$.
\end{proof}

\begin{corollary}
\label{c3.7}
For any point $P\in\mathscr P$, the divisor
$$
D^{(P)}:=\sum_{\Pi_i\ni P}\Pi_i
$$
is a cone with~vertex~$P$ over a union of four lines
forming a combinatorial cycle. In particular, the divisor~$D^{(P)}$ has four
components.
\end{corollary}

\begin{proof}
Let $H\subset X$~be a general hyperplane section. It is clear that ~$D^{(P)}$
is a cone over $H\cap D^{(P)}$ and $H\cap D^{(P)}$~is a union of
lines. By Lemma~\ref{l3.1} the divisor $H\cap D^{(P)}$ has simple
normal crossing. If $H\cap D^{(P)}$ is not connected, then~$D$ can be
decomposed in~the sum $D'+D''$ of two effective divisors so that 
$D'\cap D''=\{P\}$ (in~a neighborhood of~$P$). On the other hand, since $D$~ is a
Cartier divisor in a~variety with~terminal singularities, it is a
Cohen--Macaulay scheme that leads to a contradiction. Thus, the
intersection $H\cap D^{(P)}$ is connected.

Let $\Pi_i\subset D^{(P)}$. By Lemma~\ref{l3.6} there exist exactly two
components~$\Delta_i$ passing through~$P$. These components correspond to
two planes~$\Pi_l$, $\Pi_k$ containing~$P$. Therefore, each
component $H\cap\Pi_i\subset H\cap D^{(P)}$ intersects exactly two other
components $H\cap\Pi_l$ and $H\cap\Pi_k\subset H\cap D^{(P)}$. This 
means that $H\cap D^{(P)}$~is a combinatorial cycle.

Finally, the number of components of $H\cap D^{(P)}$ is at most $4$ by
Lemma~\ref{l3.3} and this number cannot be less than~$4$, because $H$~ is an
intersection of quadrics in~$\mathbb P^g$ and so it does not contain
``triangles'' composed of lines.
\end{proof}

\begin{lemma}
\label{l3.8}
For each plane $\Pi_i$, the intersection
$\Pi_i\cap\operatorname{Supp}(D-\Pi_i)$ has $2+a$ one-dimensional components,
where~$a$ is defined by the relation $D\sim-aK_X$ (cf.~\eqref{eq3.1}).
\end{lemma}

\begin{proof}
Let $H\subset X$~be a general hyperplane section. It is clear that $H$~ is a
nonsingular K3 surface. Let $l_i:=\Pi_i\cap H$. Since $l_i$~ is a
nonsingular rational curve, we have
$$
l_i\cdot\sum_{j\ne i}l_j=-l_i^2+l_i\cdot\sum_{j}l_j
=2+l_i\cdot\sum\Pi_i=2+a.
$$
Thus $\Pi_i$ intersects by lines exactly $2+a$ components of $D$.
\end{proof}

\begin{lemma}
\label{l3.9}
We have $|\mathscr P|=(g-1)a(a+2)(a+1)/4$.
\end{lemma}

\begin{proof}
Each plane $\Pi_i\in O$ contains $(a+2)(a+1)/2$ points from~$\mathscr P$
(which form the whole singular locus of the union of $2+a$
lines~$\Delta_i$) and there are exactly
four planes $\Pi_j\in O$ passing through each point $P\in\mathscr P$.
\end{proof}

\begin{lemma}
\label{l3.10}
We have $|\mathscr P|=\dim|D|$.
\end{lemma}

\begin{proof}
For $P\in\mathscr P$, let $H_P$~be a general hyperplane section
passing through~$P$. Let $H:=\sum_{P\in\mathscr P}H_P$ and let
$B:=(1-\delta)D+\varepsilon H$. Put
$\mathcal I_{\mathscr P}:=\mathcal I(X,B)$. For some
$0<\delta,\varepsilon\ll1$, the pair $(X,B)$ is lc and its locus
of log canonical singularities $\operatorname{LCS}(X,B)$ coincides
with~$\mathscr P$. Since the pair $(X,B)$ is lc, the scheme of
log canonical singularities is reduced. Thus
$\mathscr O_X/\mathcal I_{\mathscr P}$~is the structure sheaf of~$\mathscr P$.
Apply the Nadel Vanishing Theorem. We obtain
$H^1(X,\mathcal I_{\mathscr P}\otimes\mathscr O_X(D))=0$. Then from the exact
sequence
$$
0\longrightarrow\mathcal I_{\mathscr P}\otimes\mathscr O_X(D)
\longrightarrow\mathscr O_X(D)\longrightarrow\mathscr O_{\mathscr P}(D)
\longrightarrow0
$$
we obtain
$$
|\mathscr P|=\dim H^0(\mathscr O_{\mathscr P}(D))
=\dim H^0(\mathscr O_X(D))-\dim H^0(\mathcal I_{\mathscr P}(D)).
$$
Since $\mathscr P\subset D$, we have $H^0(\mathcal I_{\mathscr P}(D))\ne0$.
Therefore, $|\mathscr P|\le\dim|D|$.

Assume that $|\mathscr P|\le\dim|D|-1$. Let $\mathscr D\subset|D|$~ be the
linear subsystem, consisting of divisors passing through all
points of~$\mathscr P$. Then
$\dim\mathscr D\ge\dim|D|-|\mathscr P|\ge1$. Assume that planes
$\Pi_i,\Pi_j\in O$ intersect each other by a line~$l$. Then $D\cdot l=a$.
On the other hand, $l$ contains exactly $a+1$ points from~$\mathscr P$.
Therefore any element $D'\in\mathscr D$ contains all the lines of the form
$\Pi_i\cap\Pi_j$. In particular, $D'\cap\Pi_i$ contains
$\Pi_i\cap\operatorname{Supp}(D-\Pi_i)$. Since the last set~ is a 
union of $a+2$ lines, we have $D'\supset \Pi_i$ for any~$i$. Then
$D'=D$, a contradiction.
\end{proof}

\begin{proof}[Theorem \ref{t1.1}]
By the Riemann--Roch formula and Kawamata--Viehweg vanishing theorem
we have
$$
\dim|D|=\dim\,\lvert-aK_X\rvert=\frac1{12}a(a+1)(2a+1)(2g-2)+2a.
$$
Therefore, by Lemmas \ref{l3.9} and \ref{l3.10}, we obtain
$$
(a+1)(2a+1)(2g-2)+24=3(g-1)(a+2)(a+1).
$$
Thus we have
\begin{equation}
\label{eq3.2}
(a+1)(g-1)(4-a)=24.
\end{equation}
For $5\le g\le12$ the equation \eqref{eq3.2} has the following solutions:
\begin{equation}
\label{eq3.3}
(g,a,|\mathscr P|)=(5,1,6),\,(5,2,24),\,(7,3,90).
\end{equation}
The last possibility is excluded by the lemma below.
This proves our theorem.
\end{proof}

\begin{lemma}[\cite{10}]
\label{l3.11}
If $g\ge6$, then $|\operatorname{Sing}(X)|\le29$.
\end{lemma}

\begin{proof}
According to \cite[Theorem~13]{10} the number of singular points of a
Fano threefold~$X$ with~terminal singularities is at most
$$
21-\frac12\operatorname{Eu}(\mathfrak X_s)
=21-\frac12\bigl(2+ 2\operatorname{b}_2(\mathfrak X_s)
-\operatorname{b}_3(\mathfrak X_s)\bigr)
=20-\rho(\mathfrak X_s)+\operatorname{h}^{1,2}(\mathfrak X_s),
$$
where $\mathfrak X_s$~is a smoothing of $X$ as in~Theorem \ref{t2.3}. In our
case, $\rho(\mathfrak X_s)=\rho(X)=1$ and
$\operatorname{h}^{1,2}(\mathfrak X_s)\le10$ (see~\cite{9}). The lemma
is proved.
\end{proof}

For the case $\mathrm g(X)=5$, from \eqref{eq3.3} we obtain the following partial 
result.

\begin{corollary}
\label{c3.12}
Let $X$ be a $G$-Fano threefold of the main \text{series}~with $\mathrm g(X)\,{=}\,5$.
Assume that~$X$ contains a plane~$\Pi_1$ and let
$\Pi_1,\dots,\Pi_n$~be its orbit. Let $\mathscr P$~be the set of all
zero-dimensional log canonical centers of the pair $(X,\sum_i\Pi_i)$.
Then has one of the following cases holds:

{\rm(i)} $n=8$, $|\mathscr P|=6$, $|\operatorname{Sing}(X)|\ge6$;

{\rm(ii)} $n=16$, $|\mathscr P|=24$, $|\operatorname{Sing}(X)|\ge24$.
\end{corollary}

\section({\$\000\134mathbb Q\$-factorial case})
{$\mathbb Q$-factorial case}
\label{s4}

In this section we generalize the Fano--Iskovskikh ``double projection'' method to the
case of singular Fano threefolds.

\begin{theorem}
\label{t4.1}
Let $X$~be a $\mathbb Q$-factorial Gorenstein Fano threefold
of the main series with~terminal singularities and
$\mathrm g(X)\ge7$. Then there exists the following diagram:
\begin{equation}
\label{eq4.1}
\xymatrix{
& Y
\ar[dl]_f
\ar@{-->}[rr]^{\chi}
& & Y'\ar[dr]^{f'}
\\
X & & & & Z
}
\end{equation}
where $f$ is the blowup of a line $l\subset X\setminus\operatorname{Sing}(X)$,
$\chi$~is a flop, $f'$~is a Mori contraction, and
$\operatorname{Pic}(Z)\simeq\mathbb Z$.

{\rm(i)} If $g\ge9$, then $Z$~is a nonsingular Fano threefold and $f'$~ is 
the blow-up of an irreducible (possibly, singular) curve $B\subset Z$. Moreover,
\begin{table}[h]
\begin{center}
\vskip3mm
\begin{tabular}[c]{c|c| c | c}
$\mathrm g(X)$ & $Z$ & $p_a(B)$ & $-K_Z\cdot B$
\\
\hline
9 & $\mathbb P^3$ & 3 & $4\cdot7$
\\
10 & $Q\subset\mathbb P^4$\text{ -- \textit{nonsingular a quadric}} & 2 & $3\cdot7$
\\
12 & $Z_5\subset\mathbb P^6$\text{ -- \textit{nonsingular del Pezzo threefold }} & 0 & $2\cdot5$
\end{tabular}
\end{center}
\end{table}
\noindent
we have
$$
|\operatorname{Sing}(X)|=|\operatorname{Sing}(B)|\le p_a(B).
$$
In particular, $X$ is nonsingular if $\mathrm g(X)=12$.

{\rm(ii)} If $g=8$, then $f'$~is a conic bundle over
$Z\simeq\mathbb P^2$ and the discriminant curve is (possibly,
reducible) quintic $\Delta\subset\mathbb P^2$. Let $r_1$~be the number of
ordinary double points~$\Delta$, $r_2$~be the number of simple cusps, and
$r_3$~be the number of remaining singular points. Then
$$
|\operatorname{Sing}(X)|\le r_1+r_2+2r_3.
$$

{\rm(iii)} If $g=7$, then $f'$~is a del Pezzo fibration
of degree~$5$ over $Z\simeq\mathbb P^1$.
\end{theorem}

\begin{proof}
The proof  follows the classical idea of G. Fano (for a modern 
exposition for the nonsingular case we refer to~\cite{17}).

Let $X$~be a $\mathbb Q$-factorial Fano threefold
of the main series with~terminal singularities of genus
$g=\mathrm g(X)\ge7$. Then~$X$ is, in fact, factorial
\cite[Lemma~5.1]{8} and the group~$\operatorname{Cl}(X)$ is generated by the
canonical class~$K_X$. Let
$\mathfrak f\colon\mathfrak{X}\to\mathfrak{D}\ni o$~be a one-parameter
smoothing of~$X$ as in~Theorem \ref{t2.3}. By the construction, a general fiber
$X_s=\mathfrak f^{-1}(s)$~is a nonsingular Fano threefold and
$\mathfrak f^{-1}(o)=X$. According to~\cite{18} each nonsingular
fiber~$X_s$ contains a one-dimensional family of lines. Each line
deforms to a one contained in~$\mathfrak X$ and so the original variety~$X$
also contains a one-dimensional family of lines~$\mathscr L$. We claim 
that a general line~$l$ from this family~$\mathscr L$ is contained in~the nonsingular
locus of~$X$. Indeed, otherwise there exists a one-dimensional family of
lines passing through one point $P\in X$ (because the singularities of~$X$ are
isolated). All these lines swept out a surface
$F\subset X\cap\overline{T_{P,X}}$ which must be a projective cone 
over some curve. Since~$X$ is an intersection of quadrics
(by Theorem~\ref{t2.2}, (iii)) and $\dim\overline{T_{P,X}}=4$, we have
$\operatorname{deg}F\le4$. On the other hand,
$\operatorname{Cl}(X)=\mathbb Z\cdot K_X$,
a contradiction.

Thus, $X\setminus \operatorname{Sing}(X)$ contains a line~$l$.
Hereinafter, the proof goes similar to that in~\cite{17}. However, since
the threefold~$X$ can be singular in~our case, some modifications are needed. 
For convenience of the reader we present the proof completely.

Let $f\colon Y\to X$~be the blowup of $l$, $E$~be the exceptional divisor, and
let $H:=f^*(-K_X)$.

We need the following

\begin{lemma}[{\rm(see, e.g., \cite[Lemma~4.1.2]{9} or
\eqref{eq5.1})}]
\label{l4.2}
The following equalities hold:
\begin{equation}
\label{eq4.2}
(-K_{Y})^3=2g-6, \quad
(-K_{Y})^2\cdot E=3, \quad
(-K_{Y})\cdot E^2=-2, \quad
E^3=1.
\end{equation}
\end{lemma}

Using the same arguments as in \cite[Sect..~4.3.1]{9} we show 
that the linear system $\lvert-K_{Y}\rvert=|H-E|$ is nef, big, and
defines a birational morphism
$\varphi\colon Y\to Y_{2g-6}\subset\mathbb P^{g-1}$.

It is easy also to show that $\dim|H-2E|\ge g-6$ (see, e.g.,
\cite[\S\,2, Lemma~1]{17}). Since~$\operatorname{Cl}(X)$ is generated
by the class of the divisor~$-K_X$, for some $\alpha\ge2$ the linear system
$|H-\alpha E|$ has no fixed components. Using
the relations~\eqref{eq4.2}, we obtain
\begin{align*}
0\le(-K_{Y})\cdot(H-\alpha E)^2
& =(-K_{Y})\cdot(-K_{Y}-(\alpha-1)E)^2
\\
& =2g-6-6(\alpha-1)-2(\alpha-1)^2.
\end{align*}
Since $g\le12$, this gives us $\alpha=2$, i.\,e.~ the linear system
$|H-2E|$ has no fixed components.

Further we claim that $\varphi$~is a small morphism. Indeed,
otherwise~$\varphi$ contracts a prime divisor~$D$. For any fiber~$\Upsilon$
of the morphism~$\varphi$ we have $-K_{Y}\cdot\Upsilon=0$ and so
$(H-2E)\cdot\Upsilon<0$. Therefore $D$ is contained in~the base
locus of $|H-2E|$. The contradiction shows that $\varphi$~is a small
crepant morphism.

In this situation there exists a flop $\chi\colon Y\dashrightarrow Y'$, where $Y'$
has the same type of singularities as that of~$Y$ (terminal
Gorenstein)~\cite{19}. Moreover, $\rho(Y')=\rho(Y)=2$, the
divisor~$-K_{Y'}$ is nef, big, and the variety~$Y'$ (as~$Y$) 
is factorial (see~\cite[Lemma~5.1]{8}). Therefore, there exists 
an extremal Mori contraction $f'\colon Y'\to Z$. According to the general theory of
extremal rays there is the following exact sequence
\begin{equation}
\label{eq4.3}
0\longrightarrow\operatorname{Pic}(Z)\overset{{f'}^{*}}{\longrightarrow}
\operatorname{Pic}(Y')\longrightarrow\mathbb Z,
\end{equation}
where the map on the right hand side is defined by the intersection with~some curve in a~fiber.
Hence, $\rho(Z)=\rho(Y')-1=1$ and so $\dim Z>0$. Let~$H'$ and
$E'$~be the proper transforms on~$Y'$ of the divisors~$H$ and~$E$, respectively.

\begin{lemma}
\label{l4.3}
Let $F$~be a divisor on $Y'$ and $D$ be its proper transform on~$Y$.
Write $D\sim\alpha(-K_{Y})-\beta E$. Then
\begin{equation}
\begin{gathered}
\label{eq4.4}
(-K_{Y})^2\cdot D=(-K_{Y'})^2\cdot F=(2g-6)\alpha-3\beta,
\\
(-K_{Y})\cdot D^2=(-K_{Y'})\cdot{F}^2=(2g-6)\alpha^2-6\alpha\beta-2\beta^2.
\end{gathered}
\end{equation}
\end{lemma}

\begin{proof}[{\rm immediately follows from \eqref{eq4.2}}]
\end{proof}

\begin{lemma}
\label{l4.4}
Let $D$~be a prime divisor on $Y$ which is not big. For some
integers~$\alpha$ and $\beta$ we write $D\sim\alpha(-K_{Y})-\beta E$. Then
$$
\alpha,\beta>0, \quad
\beta\ge\alpha, \qquad
(-K_{Y})^2\cdot D\ge3\alpha+2\beta.
$$
\end{lemma}

\begin{proof}
Since $f(D)$ is effective, $\alpha>0$. Since the divisor~$-K_{Y}$ is big, 
$\beta>0$. The divisors~$-K_{Y'}$ and $-K_{Y'}-E'$ are nef on~$Y'$ and
they are contained in~ the closed cone of ample divisors
$\overline{\operatorname{Amp}}(Y')$. Hence~$D$ cannot be a convex
linear combination of~$-K_{Y}$ and $-K_{Y}-E$. Therefore, $\beta\ge\alpha$.
Since the linear system $\lvert-K_Y-E\rvert$ has no fixed components, we have
$$
(-K_{Y})\cdot(-K_{Y}-E)\cdot D
=(-K_{Y})^2\cdot D-(-K_{Y})\cdot E\cdot D\ge0.
$$
On the other hand, $(-K_{Y})\cdot E\cdot D=3\alpha+2\beta$
by~\eqref{eq4.2}. The lemma is proved.
\end{proof}

Below we consider the possibilities for the contraction $f'$ according to the classification 
of extremal contractions~\cite{20}.

Assume that $\dim Z=1$. Then $f'$~is a del Pezzo fibration 
and $Z\simeq\mathbb P^1$. Let $F$~be a general geometric fiber.
We use the notation of Lemma~\ref{l4.3}. Then $(-K_{Y'})\cdot F^2=0$
and $(-K_{Y'})^2\cdot F=K_{F}^2\le9$. It follows from the sequence~\eqref{eq4.3}
that $\operatorname{gcd}(\alpha,\beta)=1$ and it follows from the second relation 
in~\eqref{eq4.4} that ~$\alpha$ divides~$2$. Note that 
$-\alpha K_{F}\sim\alpha(-K_{Y'})|_{F}\sim\beta E'|_{F}$. This means that the
canonical divisor~$K_{F}$ of the del Pezzo surface~$F$ is divisible by~$\beta$.
Hence, $\beta\le3$. Moreover, if $\beta=3$, then $F\simeq\mathbb P^2$ and
$K_{F}^2=9$. In this case, taking~\eqref{eq4.4} into account we obtain $(g-3)\alpha=9$,
$\alpha=1$, and $(-K_{Y})\cdot D^2<0$. On the other hand,
$(-K_{Y})\cdot D^2=(-K_{Y'})\cdot F^2=0$, a contradiction. Similarly, if
$\beta=2$, then $K_{F}^2=8$, $(g-3)\alpha=7$, $\alpha=1$, and
$(-K_{Y})\cdot D^2<0$. Again we get a contradiction. Therefore,
$\beta=\alpha=1$ and again from~\eqref{eq4.4} we obtain $g=7$ and $K_{F}^2=5$,
i.\,e.~ the case~(iii) of our theorem.

Assume now that $\dim Z=2$. According to \cite{20} the surface~$Z$
is nonsingular and $f'$~is a conic bundle. In our case,
$\kappa(Z)=-\infty$ and $\rho(Z)=1$. Hence, $Z\simeq\mathbb P^2$. Let
$\Delta\subset\mathbb P^2$~be the discriminant curve, let
$l\subset\mathbb P^2$~be a line, and let $F:={f'}^{-1}(l)$. Again we use the
notation of Lemma~\ref{l4.3}. Since a general geometric fiber
$C\subset Y'$ is a conic, we have $(-K_{Y})\cdot C=(-K_{Y})\cdot D^2=2$ and
$0=F\cdot C=2\alpha-(E'\cdot C)\beta$. It follows from the sequence~\eqref{eq4.3}
that $\operatorname{gcd}(\alpha,\beta)=1$ and so~$\beta$
divides~$2$. By Lemma~\ref{l4.4} we have $\alpha=1$. Since
$(-K_{Y})\cdot D^2=2$, the second relation in~\eqref{eq4.4} has the form
$2g-6-6\beta-2\beta^2=2$. Hence, $\beta=1$ and $g=8$. Finally, by the
adjunction formula
$$
K_F=(K_{Y'}+F)|_F, \qquad
K_F^2= K_{Y'}^2\cdot F+2K_{Y'}\cdot F^2=3.
$$
Therefore, the projection $f'|_F\colon F\to l$ has five degenerate fibers.
Thus, $\operatorname{deg}\Delta=5$. We obtain the case~(ii) of our theorem.

Assume that the morphism $f'$ is birational and contracts an (irreducible)
divisor~$F$ to a~point. Let, as above, $D\subset Y$~be the proper
transform~$F$ and $D\sim\alpha(-K_{Y})-\beta E$. According to the classification
from~\cite{20} there exist four types of such contractions and in all these cases
$(-K_{Y'})^2\cdot D'\le4$. This contradicts Lemma~\ref{l4.4}.

Finally, assume that the morphism $f'$ is birational and contracts an
(irreducible) divisor~$F$ to a curve~$B$. According to~\cite{20} the singularities of the
curve~$B$ are locally planar, the variety~$Z$ is nonsingular along~$B$, and~$f'$
is the blowup of the ideal sheaf of~$B$. Then $Z$~is a Fano threefold
with~terminal factorial singularities and $\rho(Z)=1$. Let $A$~ be the
positive generator of the group~$\operatorname{Pic}(Z)$. Then
$-K_Z=\iota A$ for some positive integer~$\iota$ which
is called the \textit{Fano index}. It is well-known that $\iota\le4$
(see~\cite{9} and Theorem~\ref{t2.3}). Moreover, $\iota=4$ if and
only if $Z\simeq\mathbb P^3$, and $\iota=3$ if and only
if $Z$~is a quadric in~$\mathbb P^4$.

Below we use the notation of Lemma \ref{l4.3}. Let $C$~be a general
fiber of $f'|_F\colon F\to B$. Since over a general point of the curve~$B$ the morphism
$f'$~is a usual blowup, $F\cdot C=-1$. Therefore,
$(E'\cdot C)\beta=\alpha+1$. In particular, $\beta$ divides $\alpha+1$. Since
$\dim|F|=0$ and $\dim\lvert-K_{Y'}-E'\rvert>0$, we have $\alpha\ne\beta$. Then from
Lemma~\ref{l4.4} we obtain $\beta>\alpha$. Hence, $\beta=\alpha+1$. Further,
\begin{gather*}
K_{Y'}=(f')^*K_Z+F=-\iota(f')^*A+\alpha(-K_{Y'})-(\alpha+1)E',
\\
\iota{f'}^*A=(\alpha+1)(-K_{Y'}-E').
\end{gather*}
Since the divisors $(f')^*A$ and $-K_{Y'}-E$~are primitive elements of the lattice
$\operatorname{Pic}(Y')$, we have $\beta=\alpha+1=\iota$ and
$(f')^*A=-K_{Y'}-E'$. In particular, $1\le\alpha\le\iota-1=3$. Moreover,
\begin{equation}
\label{eq4.5}
\dim|A|\ge\lvert-K_{Y}-E\rvert\ge g-6.
\end{equation}
The intersection theory on $Y'$ has the same form as the intersection theory on 
the blowup of a nonsingular variety along a nonsingular curve
(see~\eqref{eq5.1}). Hence,
\begin{gather*}
(-K_{Y})^2\cdot D=-K_Z\cdot B-2p_a(B)+2,
\\
(-K_{Y})\cdot D^2=2p_a(B)-2.
\end{gather*}
Taking \eqref{eq4.4} and $\beta=\alpha+1$ into account we obtain
\begin{gather*}
(2g-6)\alpha-3(\alpha+1)=-K_Z\cdot B-2p_a(B)+2,
\\
(2g-6)\alpha^2-6\alpha(\alpha+1)-2(\alpha+1)^2=2p_a(B)-2.
\end{gather*}
Adding up the last two equalities we obtain
$$
(2g-6)\alpha-3(\alpha+1)+(2g-6)\alpha^2-6\alpha(\alpha+1)-2(\alpha+1)^2
=-K_Z\cdot B.
$$
Since
$\beta=\alpha+1=\iota$, we have
$$
2(g-7)\alpha=A\cdot B+5.
$$
Consider cases $\alpha=1,2,3$ separately.

Let $\alpha=1$. Then $\beta=\iota=2$ and
\begin{equation}
\label{eq4.6}
(2g-6)\alpha^2-6\alpha\beta-2\beta^2=2g-26=2p_a(B)-2.
\end{equation}
The only solution for \eqref{eq4.6} is $p_a(B)=0$, $g=12$,
$A\cdot B=5$. In this case, ~$Z$ is a del Pezzo threefold (see,
e.g., \cite{9} or~\cite{2}). According to~\eqref{eq4.5} we have
$\dim|A|\ge6$. Since $\rho(Z)=1$ and~$Z$ factorial, $A^3=5$.
According to \cite[Corollary~5.4]{2} the variety~$Z$ is nonsingular.

Let $\alpha=2$. Then $\beta=\iota=3$ and $Z$~is a quadric in~$\mathbb P^4$.
Since the variety~$Z$ is factorial, this quadric is nonsingular. As above,
we have
\begin{equation}
\label{eq4.7}
(2g-6)\alpha^2-6\alpha\beta-2\beta^2=8g-78=2p_a(B)-2, \quad
4g=p_a(B)+38.
\end{equation}
According to \eqref{eq4.5} we have $4=\dim|A|\ge g-6$. Hence, $g\le10$.
The only solution of \eqref{eq4.7} is 
$g=10$, $p_a(B)=2$, $A\cdot B=7$.

Let $\alpha=3$. Then $\beta=\iota=3$ and $Z\simeq\mathbb P^3$. As above,
$9g=p_a(B)+78$, $3=\dim|A|\ge g-6$, $g=9$, $p_a(B)=3$, $A\cdot B=7$.

Thus the existence of the diagram \eqref{eq4.1} and its properties are proved.

For the proofs of the assertion about singularities we have to notice that, by the
construction, ~$f$ is an isomorphism near~$\operatorname{Sing}(X)$
and~$\operatorname{Sing}(Y)$, and the map~$\chi$ preserves completely the type of
singularities (and their number)~\cite{19}. Thus,
$|\operatorname{Sing}(X)|=|\operatorname{Sing}(Y')|$. The bound for
$|\operatorname{Sing}(Y')|$ follows from Proposition~\ref{p5.2} in~ the case
$g=8$ and from Proposition~\ref{p5.1} in cases $g\ge9$. The theorem is proved.
\end{proof}

\begin{remark}
\label{r4.5}
For $g\ge9$, the construction \eqref{eq4.1} can be reversed: for a suitable
choice of the curve~$B$ with~corresponding values of degree and arithmetic
genus its the blowup $f'\colon Y'\to Z$ satisfies the standard
conditions:

a) the linear system $\lvert-K_Y'\rvert$ is base point free;

b) the corresponding morphism $\Phi_{\lvert-K_Y'\rvert}$ does not contract any divisors;

c) the variety $Y'$ has only terminal singularities.

In this situation, there exists (and reconstructed uniquely) the right hand part of the 
diagram~\eqref{eq4.1}. Such a curve can be chosen on a nonsingular
del Pezzo surface of degree~$3$, $4$, $5$ in~cases $g=9,10,12$,
respectively. This allows to resolve problems on the existence of
Fano threefolds with~ given number of singular points.

In particular, this construction allows to construct examples of non-projective 
Moishezon threefolds with~$b_2=1$ as small resolutions of
our variety~$X$ (cf.~\cite{21}).
\end{remark}

\begin{proof}[Theorem \ref{t1.3}]
In the cases $g\ge9$ the assertion immediately follows from
Theorem~\ref{t4.1}, (i). Consider the case $g=8$. For a plane reduced
(but possibly reducible) curve~$C$ put $\gamma(C):=r_1+r_2+2r$,
where~$r_1$ (respectively, $r_2$)~is the number of simple double points of type~$A_1$
(respectively, the number of double points of type $A_2$), 
and $r$~is the number of remaining singular points. Then
by Theorem~\ref{t4.1}, (ii) we have $|\operatorname{Sing}(X)|\le\gamma(\Delta)$.
The estimate $\gamma(\Delta)\le10$ for a plane quintic~$\Delta$ follows from
the following two simple assertions:

1) if the curve $C$ is irreducible, then $\gamma(C)\le p_a(C)$;

2) if $C_1$~is an irreducible nonsingular component of $C$, then
$\gamma(C)\le\gamma(C-C_1)+C\cdot(C-C_1)$.

The theorem is proved.
\end{proof}

\begin{remark}
\label{r4.6}
In contrast with the case $g\ge9$, we do not assert that the bound
$|\operatorname{Sing}(X)|\le10$ is sharp for $g=8$. One can conjecture that 
it can be improved.
\end{remark}

\section{Two auxiliary results}
\label{s5}

\begin{prop}
\label{p5.1}
Let $V$~be a threefold with~terminal
singularities and let $f\colon V\to W$~be a birational Mori contraction that
contracts a divisor~$F$ to a curve~$B$. Then the following
assertions hold:

{\rm(i)} the singularities of the curve $B$ are locally planar, the variety~$W$ is nonsingular
along~$B$, and~$f$ is the blowup of the ideal sheaf of~$B$;

{\rm(ii)} $f(\operatorname{Sing}(X)\cap F)=\operatorname{Sing}(B)$ and
each fiber $f^{-1}(b)$ over a point $b\in\operatorname{Sing}(B)$ contains
exactly one singularity of~$X$.

Moreover,
\begin{equation}
\begin{gathered}
\label{eq5.1}
(-K_V)^3=(-K_W)^3+2K_W\cdot B+2p_a(B)-2,
\\
(-K_V)^2\cdot E=-K_W\cdot B-2p_a(B)+2,
\\
(-K_V)\cdot E^2=2p_a(B)-2.
\end{gathered}
\end{equation}
\end{prop}

\begin{proof}
The assertion (i) follows from \cite{20} and the assertion (ii)~is 
a simple computation in~local coordinates. Let us prove~\eqref{eq5.1}. For some
ample divisor~$A$ on~$W$, the linear system $\lvert-K_V+f^*A\rvert$ has no 
base points (see~\cite[Proposition~1]{20}). Take a general element
$S\in\lvert-K_V+f^*A\rvert$. By Bertini's theorem~$S$ is nonsingular. 
Let $\bar S:=f(S)$.
Then $\bar S\in\lvert-K_W+A\rvert$ and $f^*\bar S=S+E$. The restricted linear
system $\lvert-K_V+f^*A\rvert\bigr|_E$ is ample and base point free. Again by
Bertini's theorem its general element $S\cap E$~is a nonsingular irreducible curve.
Since the intersection number of~$S$ and a~general fiber $E\to f(E)$ equals ~$1$,
the restriction $f_S\colon S\to\bar S$ is an isomorphism and
$f_S(E\cap S)=B$. Note that $K_S=f^*A|_S$, $K_{\bar S}=A|_{\bar S}$, and
$(B\cdot B)_{\bar S}=2p_a(B)-2-A\cdot B$. Using the last relation we can
write
\begin{gather*}
-K_V\cdot E^2=(S-f^*A)\cdot E^2=(B\cdot B)_{\bar S}+A\cdot B=2p_a(B)-2,
\\
-K_V\cdot f^*K_W\cdot E=(S-f^*A)\cdot f^*K_W\cdot E
=S\cdot f^*K_W\cdot E=K_W\cdot B.
\end{gather*}
Hence we have
\begin{gather*}
K_V^2\cdot E=K_V\cdot f^*K_W\cdot E+K_V\cdot E^2=-K_W\cdot B-2p_a(B)+2,
\\
\begin{split}
(-K_V)^3
&=-K_V\cdot(f^*K_W+E)^2=-(f^*K_W+E)\cdot f^*K_W^2
\\
&\qquad
-2K_V\cdot f^*K_W\cdot E-K_V\cdot E^2=(-K_W)^3+2K_W \cdot B+2p_a(B)-2.
\end{split}
\end{gather*}
The proposition is proved.
\end{proof}

\begin{prop}
\label{p5.2}
Let $V$~be a threefold with~terminal
singularities and let $f\colon V\to W$~be a Mori contraction to a surface.
Then the surface~$W$ is nonsingular and $f$~is a conic bundle
(possibly, singular). Further, let $\Delta\subset W$~be the discriminant
curve. Then
$f(\operatorname{Sing}(V))\subset\operatorname{Sing}(\Delta)$. Moreover,
any fiber $f^{-1}(w)$, $w\in\operatorname{Sing}(\Delta)$ contains at most two
points of $\operatorname{Sing}(V)$. If
$f^{-1}(w)\cap\operatorname{Sing}(V)$~consists of exactly two points, then the singularity
$w\in\Delta$ is not an ordinary double point~$A_1$ nor a simple cusp~$A_2$.
\end{prop}

\begin{proof}
The first part of the proposition is contained in~\cite{20}. 
It remains to prove only the assertion about singularities of~$V$. Since 
the problem is local, we may
assume that the ground field ~$\Bbbk$~is the field of complex numbers~$\mathbb C$,
$V$~is an analytic neighborhood of a fiber~$f^{-1}(w)$, and
$W\subset\mathbb C^2_{u,v}$~is a small disk containing $w=(0,0)$. Then
we can embed~$V$ to~$\mathbb P^2\times W$ so that ~$V$ is defined
by the equation $q(x,y,z;u,v)=0$, where~$q$ is regarded as a quadratic form in
~$x$, $y$, $z$ with~coefficients in~$\mathbb C\{u,v\}$. The fiber~$f^{-1}(w)$
is defined by the equation $q(x,y,z;0,0)=0$. Since $f^{-1}(w)$~is a conic, we have
$\operatorname{rk}q(x,y,z;0,0)\ge1$. If
$\operatorname{rk}q(x,y,z;0,0)=3$, then the fiber~$f^{-1}(w)$ is nonsingular and~$V$ is also
nonsingular (near~$f^{-1}(w)$). If $\operatorname{rk}q(x,y,z;0,0)=2$, then
up to coordinate change we can write $q(x,y,z;0,0)=x^2+y^2$ and
$q(x,y,z;u,v)=x^2+y^2+\alpha(u,v)z^2$, where $\alpha=0$~is the equation of~$\Delta$
and $\alpha(0,0)=0$. In this case,~$V$ is singular, if and only if
$\operatorname{mult}_{(0,0)}\alpha>1$, i.\,e. the~curve~$\Delta$ is singular at~the origin. Moreover,
$\operatorname{Sing}(V)\subset\operatorname{Sing}(f^{-1}(w))=\{P\}$, where
$\operatorname{Sing}(f^{-1}(w))$~is a single point.

Finally, consider the case $\operatorname{rk}q(x,y,z;0,0)=1$. Then
up to coordinate change we can write $q(x,y,z;0,0)=x^2$ and
$q(x,y,z;u,v)=x^2+\alpha y^2+2\beta yz+\gamma z^2$,
where~$\alpha$, $\beta$, $\gamma$~are holomorphic functions in~$u$, $v$,
vanishing at~the origin. The equation of ~$\Delta$ has the form
$\alpha\gamma-\beta^2=0$. Hence the curve~$\Delta$ is singular at~$(0,0)$.
Assume that ~$V$ has two singular points~$P_1$, $P_2$ on~$f^{-1}(w)$.
By changing the coordinates~$y$, $z$ linearly, we may assume that $P_1=(0:1:0)$,
$P_2=(0:0:1)$. Then $\operatorname{mult}_{(0,0)}\alpha>1$ and
$\operatorname{mult}_{(0,0)}\gamma>1$. Since the singularities of~$V$
are isolated, we have $\operatorname{mult}_{(0,0)}\beta=1$. Then it is easy to see
that the variety~$V$ is nonsingular outside of ~$P_1$, $P_2$ and
the singularity
$\{\alpha\gamma-\beta^2=0\}$ 
is not an ordinary double point nor a simple cusp.
\end{proof}

\end{fulltext}

\end{document}